\documentclass[3p,times]{amsart}

\usepackage{hyperref}
\usepackage{amsmath, amssymb}
\usepackage{amscd}
\usepackage{xcolor}
\usepackage{amsthm}
\usepackage[font=footnotesize]{caption}

\usepackage{enumerate} 

\newtheorem{theorem}{Theorem}[section] 
\newtheorem{lemma}[theorem]{Lemma}   
\newtheorem{corollary}[theorem]{Corollary}
\newtheorem{proposition}[theorem]{Proposition}
\newtheorem{definition}[theorem]{Definition}

\newtheorem{main-theorem}[theorem]{Theorem}

\newtheorem*{problem*}{Problem}

\theoremstyle{definition}

\newtheorem*{question*}{Question}

\newtheorem{remark}[theorem]{Remark}

\newcommand{\wt}{\widetilde}

\newcommand{\cA}{\mathcal{A}}
\newcommand{\cU}{\mathcal{U}}

\newcommand{\cX}{\mathcal{X}}
\newcommand{\La}{\Lambda}
\newcommand{\ba}{\bar{\alpha}}

\newcommand{\ve}{\varepsilon}
\newcommand{\cV}{\mathcal{V}}

\makeatletter
\newcommand{\vect}[1]{%
  \vbox{\m@th \ialign {##\crcr
  \vectfill\crcr\noalign{\kern-\p@ \nointerlineskip}
  $\hfil\displaystyle{#1}\hfil$\crcr}}}
\def\vectfill{%
  $\m@th\smash-\mkern-7mu%
  \cleaders\hbox{$\mkern-2mu\smash-\mkern-2mu$}\hfill
  \mkern-7mu\raisebox{-3.81pt}[\p@][\p@]{$\mathord\mathchar"017E$}$}

\newcommand{\amsvect}{%
  \mathpalette {\overarrow@\vectfill@}}
\def\vectfill@{\arrowfill@\relbar\relbar{\raisebox{-3.81pt}[\p@][\p@]{$\mathord\mathchar"017E$}}}

\newcommand{\amsvectb}{%
  \mathpalette {\overarrow@\vectfillb@}}
\newcommand{\vecbar}{%
  \scalebox{0.8}{$\relbar$}}
\def\vectfillb@{\arrowfill@\vecbar\vecbar{\raisebox{-4pt}[\p@][\p@]{$\mathord\mathchar"017E$}}}
\makeatother
\usepackage{etex} 

\usepackage{etex} 

\usepackage{dsfont}         
\usepackage[h]{esvect}

\usepackage[matrix,arrow,curve,frame,arc]{xy}

\usepackage{pgf}
\usepackage{tikz}

\usepackage{color}

\usetikzlibrary{arrows,shapes,automata,backgrounds,petri,calc}

\newcommand{\tikzAngleOfLine}{\tikz@AngleOfLine}
\def\tikz@AngleOfLine(#1)(#2)#3{%
\pgfmathanglebetweenpoints{%
\pgfpointanchor{#1}{center}}{%
\pgfpointanchor{#2}{center}}
\pgfmathsetmacro{#3}{\pgfmathresult}%
}

\begin{document}

\baselineskip=14pt

\title{Some cluster tilting modules for weighted surface algebras}

\author[K. Erdmann]{Karin Erdmann}
\address[Karin Erdmann]{Mathematical Institute,
   Oxford University,
   ROQ, Oxford OX2 6GG,
   United Kingdom}
\email{erdmann@maths.ox.ac.uk}

\begin{abstract}
Non-singular weighted surface algebras satisfy the necessary condition found in  \cite{EH} 
for existence of cluster tilting modules.  We show that  any such algebra whose Gabriel quiver is bipartite,
has a module satisfying the necessary ext vanishing condition. We show that it is
3-cluster tilting precisely for non-singular triangular or spherical algebras, but not for any other
weighted surface algebra with bipartite  Gabriel quiver.

\bigskip



\end{abstract}


\maketitle

\bigskip

\section{Introduction}

A module $M$ of a finite-dimensional algebra $A$ is an $n$-cluster tilting module (or maximal $(n-1)$-orthogonal) provided
\begin{align*} {\rm add}(M)=& \{ N\mid {\rm Ext}^i(M, N)=0 \mbox{ for } 1\leq i\leq n-1\}\cr
              = &\{ N\mid {\rm Ext}^i(N, M)=0 \mbox{ for } 1\leq i\leq n-1\},
\end{align*}             
(see \cite{I1}, \cite{I2}).
We would like to know whether non-singular weighted surface algebras  have cluster tilting modules.
Weighted surface algebras  are a class of tame symmetric algebras, periodic as bimodules, of period $4$ (see \cite{WSA} and \cite{WSA-GV}, \cite{WSA-corr}).
This means that they satisfy the necessary condition found in \cite{EH}, requiring that all non-projective modules should have bounded periodic resolutions. 
As observed in \cite{EH} if such an algebra  has an $n$-cluster tilting module then  the only option is $n=3$.

Here we study weighted surface algebras which have a bipartite Gabriel quiver, which means that in the presentation as in \cite{WSA-GV} (see also  \cite{WSA-corr}) 
it has many virtual arrows.  We introduce  a module $M$, defined in
the same way for each of the algebras, which satisfies ${\rm Ext}^1(M, M)=0$ and ${\rm Ext}^2(M, M)=0$. We show that it is 3-cluster tilting 
when $\La$ is either a triangle algebra $T(\lambda)$, or a spherical algebra $S(\lambda)$ (see \S 3 for the definition) with $\lambda \in K$ and 
$\lambda\neq 0, 1$. We also show that for any other 
weighted surface algebra whose Gabriel quiver is bipartite, $M$ cannot be a direct summand of a 3-cluster tilting module.

\medskip

The algebra $T(\lambda)$ occurs in various places in the literature.  It is an algebra with $k=1$ 
in the family $Q(3\cA)_1$ of algebras of quaternion type, in  \cite{E}.  Furthermore,  it occurs with the name  $B_{1,1}(\lambda)$ in \cite{BIKR}. 
 As well, it occurs in \cite{BS}  with the name $A_1(\lambda)$. 
 In \cite{WSA-GV} it is called  the triangle algebra $T(\lambda)$ (in Example 3.4).
 Similarly the spherical algebra $S(\lambda)$ was introduced in \cite{WSA-GV}
 (Example 3.6).

Spherical algebras 
 are a special case of the family of algebras which come from the
 triangulation $T(n)$ of the sphere as defined in Example 7.5 of \cite{AGDT}. We call these algebras $n$-spherical; when $n=2$ they
 are the same as the spherical algebras. One may also observe that the spherical algebra with $n=1$ 
with the multiplicities $2, 2, 1$  is  the same as the triangle algebra (see \cite{WSA-GV}, Example 3.3).

 Algebras whose  Gabriel quiver is the same as that of $T(\lambda)$, allowing 
  a multiplicity $k>1$, coincide (up to 
 a scalar parameter) with the algebras $Q(3\cA)_2^k$ in the labelling of \cite{H}. When the characteristic of $K$ is $2$ these occur
 as the basic algebras of blocks of finite groups. 
  Very recently B. B\"{o}hmler and R. Marcinczik  proved using computer calculations that  for $k=2$, it has a  3-cluster tilting module (see \cite{BM}).

 Much of this note was written five years ago,  when  talking to  Idun Reiten
 about \cite{BIKR}, and it was
 extended first when spherical algebras had been discovered, and then again, inspired by an email from R. Marcinczik (for which I am  grateful).

 \section{Preliminaries}
 
  Throughout $K$ is an algebraically closed field, of arbitrary characteristic.
 Assume $\La$ is a finite-dimensional symmetric $K$-algebra.  
 We recall some identities for the stable category $\underline{\rm mod}\La$.

(1) \ $D{\rm Ext}^1(M, N)\cong \underline{\rm Hom}(\tau^{-1}N, M)$, and in this case $\tau\cong \Omega^2$.

(2) \ ${\rm Ext}^i(U, V)\cong \underline{\rm Hom}(\Omega^iU, V)\cong \underline{\rm Hom}(U, \Omega^{-i}V)$. 

This implies that  $\dim {\rm Ext}^i(M, N) = \dim {\rm Ext}^1(N, \Omega^{i+2}M)$. 
The algebras we consider  have the property that
all non-projective indecomposable right $\La$-modules are $\Omega$ periodic of periods dividing $4$. 
This gives us the following, we refer to this as ext symmetry.

\bigskip

\begin{corollary}  Assume $\La$ is symmetric and all modules have $\Omega$-period dividing $4$. Then 
for all $M, N$ we have $\dim {\rm Ext}^2(M, N)= \dim {\rm Ext}^1(N, M)$ as vector spaces. 
\end{corollary}

This simplifies the search for 3-cluster tilting modules.
If we know that 
${\rm Ext}^1(N,X)=0$ and ${\rm Ext}^1(X, N)=0$ then automatically ${\rm Ext}^2(X, N)=0$ and 
${\rm Ext}^2(N,X)=0$.

\medskip

\section{The algebras}

\subsection{Weighted surface algebras}

We review the definition from \cite{WSA-corr}, for details see \cite{WSA}, \cite{WSA-GV}, \cite{WSA-corr}.

Assume  $Q$ is a finite quiver.
Denote by $KQ$ the path algebra of $Q$ over $K$.
We will consider algebras of the form $A=KQ/I$ where $I$ is an ideal of $KQ$ 
which contains all paths of length $\geq m$ for some $m>>0$, so that 
the algebra is finite-dimensional and basic.
The Gabriel quiver $Q_A$ of $A$ is then the full subquiver of $Q$ 
obtained from $Q$ by removing all arrows $\alpha$ 
with $\alpha + I \in R_Q^2 + I$.

A quiver $Q$ is \emph{$2$-regular} if for each vertex $i\in Q_0$ 
there are precisely two arrows starting
at $i$ and two arrows ending at $i$. 
Such a quiver has an involution
on the arrows, $\alpha \mapsto \ba$,  
such that for each arrow $\alpha$, 
the arrow $\ba$ is the arrow $\neq \alpha$ such that $s(\alpha) = s(\ba)$. 
 
  A \emph{triangulation   quiver} \  is a pair $(Q, f)$ where $Q$ 
is a  (finite) connected 2-regular quiver,  with at least two vertices, and where 
$f$ is a fixed 
permutation of the arrows such that  $t(\alpha) = s(f(\alpha))$ 
for each arrow $\alpha$, and such that $f^3$ is the identity. 
The permutation $f$ uniquely determines a permutation $g$ of the arrows, 
defined by $g(\alpha) := \overline{f(\alpha)}$ for any arrow $\alpha$. 
We assume throughout that $(Q, f)$ is a triangulation quiver. To give the presentations of the algebras in question, we
use the following notation. 
For each arrow $\alpha$, we fix 
\begin{align*} 
m_{\alpha}\in \mathbb{N}^* &&& \mbox{ a weight, constant on $g$-cycles, and }\cr
c_{\alpha} \in K^*  &&&  \mbox{ a parameter,  constant on $g$-cycles, and define }\cr
n_{\alpha}:=   &&&  \mbox{ the length of the $g$-cycle of $\alpha$, } \cr
B_{\alpha}:=  \alpha g(\alpha)\ldots g^{m_{\alpha}n_{\alpha}-1}(\alpha) &&&   \mbox{ the path  along the $g$-cycle of $\alpha$ 
                     of length $m_{\alpha}n_{\alpha}$}, \cr
 A_{\alpha}:=  \alpha g(\alpha)\ldots g^{m_{\alpha}n_{\alpha}-2}(\alpha) &&&   \mbox{ the path  along the $g$-cycle of $\alpha$ 
                     of length $m_{\alpha}n_{\alpha}-1$.}
\end{align*}

\begin{definition}\label{def:virtual} \normalfont We say that an arrow $\alpha$ of
	$Q$ is virtual if $m_{\alpha}n_{\alpha}=2$, that is $A_{\alpha}$ has length $1$. Note that
	this condition is preserved under the permutation $g$, and that
	virtual arrows form $g$-orbits of sizes 1 or 2. 
\end{definition}


We  assume that 
the following conditions hold.\\
(1) \ 
	$m_{\alpha}n_{\alpha}\geq 2$ for all arrows $\alpha$, and\\
(2) \ 	$m_{\alpha}n_{\alpha}\geq 3$ for all arrows $\alpha$ such
	that $\ba$ is virtual and  $\ba$ is not a loop, and
	$m_{\alpha}n_{\alpha}\geq 4$ for all arrows $\alpha$ 
	such that $\ba$ is virtual and  $\ba$ is a loop.\\
Condition (1) is a general assumption, and (2) is needed to eliminate
	two small algebras (see \cite{WSA-GV}).  
We also assume that $Q$ has at least three vertices. 
With this, the definition of a weighted surface algebra (as revised in 
\cite{WSA-corr}) is as follows.

\begin{definition} \label{def:2.2}
The algebra $\La = \La(Q, f, m_{\bullet}, c_{\bullet}) = KQ/I$ 
is a weighted surface algebra if
$(Q, f)$ is a triangulation quiver, with $|Q_0| \geq 2$, 
and $I= I(Q, f, m_{\bullet}, c_{\bullet})$ is the ideal of 
$KQ$ generated by:
\begin{enumerate}[\rm(1)]
 \item
	$\alpha f(\alpha) - c_{\ba}A_{\ba}$ for all arrows
	$\alpha$ of $Q$,
 \item
	$\alpha f(\alpha) g(f(\alpha))$ \ for all arrows $\alpha$ of $Q$
	unless  $f^2(\alpha)$ is virtual, or unless $f(\ba)$ is virtual and $m_{\ba}=1, \  n_{\ba}=3$. 
 \item
	$\alpha g(\alpha)f(g(\alpha))$ for all arrows $\alpha$ of $Q$
	unless $f(\alpha)$ is virtual, or unless $f^2(\alpha)$ is virtual and $m_{f(\alpha)}=1, \  n_{f(\alpha)}=3$. 
\end{enumerate}
\end{definition}

The Gabriel quiver $Q_{\La}$ is the subquiver of $Q$ obtained by removing all
virtual arrows.\\

We recall a few  properties. \\
(1) Any such algebra is symmetric and tame.\\
(2) The dimension of $e_i\La$ is equal to $m_{\alpha}n_{\alpha} + m_{\ba}n_{\ba}$ where $\alpha, \ba$ are
the arrows starting at $i$.\\
The relations also imply that $c_{\alpha}B_{\alpha} = c_{\ba}B_{\ba}$ in $\La$. One can show that this  spans the socle of $e_i\La$.

\bigskip

We wish to define a module $M$ such that ${\rm Ext}^1(M, M)=0$ and ${\rm Ext}^2(M, M)=0$, as a candidate to be 3-cluster tilting. 
This can be done for a weighted surface algebra
whose quiver is bipartite; this   requires that each triangle of $f$ must contain a virtual arrow.
Such a quiver can be thought of made up of three building blocks, first 
a quiver of the form
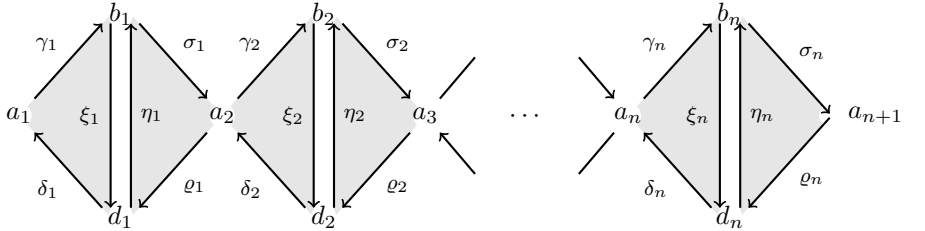
\[
\begin{tikzpicture}
[->,scale=.9]
\coordinate (0) at (-1.5,0);
\coordinate (1) at (0,1.5);
\coordinate (1l) at (-0.15,1.5);
\coordinate (1p) at (0.15,1.5);
\coordinate (2) at (0,-1.5);
\coordinate (2l) at (-0.15,-1.5);
\coordinate (2p) at (0.15,-1.5);
\coordinate (3) at (1.5,0);
\coordinate (4) at (3,1.5);
\coordinate (4l) at (2.85,1.5);
\coordinate (4p) at (3.15,1.5);
\coordinate (5) at (3,-1.5);
\coordinate (5l) at (2.85,-1.5);
\coordinate (5p) at (3.15,-1.5);
\coordinate (6) at (4.5,0);
\coordinate (7) at (6,1.5);
\coordinate (7l) at (5.85,1.5);
\coordinate (7p) at (6.15,1.5);
\coordinate (8) at (6,-1.5);
\coordinate (8l) at (5.85,-1.5);
\coordinate (8p) at (6.15,-1.5);
\coordinate (9) at (7.5,0);
\coordinate (10) at (9,1.5);
\coordinate (10l) at (8.85,1.5);
\coordinate (10p) at (9.15,1.5);
\coordinate (11) at (9,-1.5);
\coordinate (11l) at (8.85,-1.5);
\coordinate (11p) at (9.15,-1.5);
\coordinate (12) at (10.5,0);
\fill[fill=gray!20] (1l) -- (2l) -- (0) -- cycle;
\fill[fill=gray!20] (1p) -- (2p) -- (3) -- cycle;
\fill[fill=gray!20] (4l) -- (5l) -- (3) -- cycle;
\fill[fill=gray!20] (4p) -- (5p) -- (6) -- cycle;
\fill[fill=gray!20] (10l) -- (11l) -- (9) -- cycle;
\fill[fill=gray!20] (10p) -- (11p) -- (12) -- cycle;
\node[fill=white,circle,minimum size=3]  (0) at (-1.5,0) { };
\node (0) at (-1.5,0) {$a_1$};
\node[fill=white,circle,minimum size=3]  (1) at (0,1.5) { };
\node (1) at (0,1.5) {$b_1$};
\node (1l) at (-0.15,1.5) { };
\node (1p) at (0.15,1.5) { };
\node[fill=white,circle,minimum size=3]  (2) at (0,-1.5) { };
\node (2) at (0,-1.5) {$d_1$};
\node (2l) at (-0.15,-1.5) { };
\node (2p) at (0.15,-1.5) { };
\node[fill=white,circle,minimum size=3]  (3) at (1.5,0) { };
\node (3) at (1.5,0) {$a_2$};

\node[fill=white,circle,minimum size=3]  (4) at (3,1.5) { };
\node (4) at (3,1.5) {$b_2$};
\node (4l) at (2.85,1.5) { };
\node (4p) at (3.15,1.5) { };
\node[fill=white,circle,minimum size=3]  (5) at (3,-1.5) { };
\node (5) at (3,-1.5) {$d_2$};
\node (5l) at (2.85,-1.5) { };
\node (5p) at (3.15,-1.5) { };
\node[fill=white,circle,minimum size=3]  (6) at (4.5,0) { };
\node (6) at (4.5,0) {$a_3$};

\node (7l) at (5.35,1.) { };
\node (7p) at (6.65,1.) { };
\node (8l) at (5.35,-1.) { };
\node (8p) at (6.65,-1.) { };
\node at (6,0) {$\cdots$};

\node[fill=white,circle,minimum size=3]  (9) at (7.5,0) { };
\node (9) at (7.5,0) {$a_n$};
\node[fill=white,circle,minimum size=3]  (10) at (9,1.5) { };
\node (10) at (9,1.5) {$b_n$};
\node (10l) at (8.85,1.5) { };
\node (10p) at (9.15,1.5) { };
\node[fill=white,circle,minimum size=3]  (11) at (9,-1.5) { };
\node (11) at (9,-1.5) {$d_n$};
\node (11l) at (8.85,-1.5) { };
\node (11p) at (9.15,-1.5) { };
\node[fill=white,circle,minimum size=3]  (12) at (10.5,0) { };
\node (12) at (10.5,0) {$a_{n+1} \!\!\!\!\!\!\!\!\!\!\!\!\!\!\!\!\!\!\!\!$};

\fill[fill=gray!20] (1l) -- (2l) -- (0) -- cycle;
\fill[fill=gray!20] (1p) -- (2p) -- (3) -- cycle;

\draw[thick,->]
(1l) edge node[left]{\footnotesize$\xi_1$} (2l)
(2l) edge node[below left]{\footnotesize$\delta_1$} (0)
(0) edge node[above left]{\footnotesize$\gamma_1$} (1l)
(2p) edge node[right]{\footnotesize$\eta_1$} (1p)
(1p) edge node[above right]{\footnotesize$\sigma_1$} (3)
(3) edge node[below right]{\footnotesize$\varrho_1$} (2p)
(4l) edge node[left]{\footnotesize$\xi_2$} (5l)
(5l) edge node[below left]{\footnotesize$\delta_2$} (3)
(3) edge node[above left]{\footnotesize$\gamma_2$} (4l)
(5p) edge node[right]{\footnotesize$\eta_2$} (4p)
(4p) edge node[above right]{\footnotesize$\sigma_2$} (6)
(6) edge node[below right]{\footnotesize$\varrho_2$} (5p)
(8l) edge (6)
(7p) edge (9)
(10l) edge node[left]{\footnotesize$\xi_n$} (11l)
(11l) edge node[below left]{\footnotesize$\delta_n$} (9)
(9) edge node[above left]{\footnotesize$\gamma_n$} (10l)
(11p) edge node[right]{\footnotesize$\eta_n$} (10p)
(10p) edge node[above right]{\footnotesize$\sigma_n$} (12)
(12) edge node[below right]{\footnotesize$\varrho_n$} (11p)
;
\draw[thick,-]
(6) edge (7l)
(9) edge (8p)
;
\end{tikzpicture}
\qquad \quad \raisebox{9ex}{,}
\]
where the shaded triangles define the  $f$-orbits.

Next,  quivers of the form

\[
 \xymatrix{
  1 \ar@(dl,ul)[]^{\ve} \ar@<+.5ex>[r]^{\alpha}
   & 2 \ar@<+.5ex>[l]^{\beta} 
 }
\]
or 
\[
 \xymatrix{
  2'
    \ar@<+.5ex>[r]^{\gamma}
   & 3  \ar@<+.5ex>[l]^{\delta} \ar@(ur,dr)[]^{\ve'}
 }
\]

We describe the quivers and algebras we consider.
We always take the multiplicities at $2$-cycles of $g$ equal to $1$, and at loops we take multiplicity  $2$. 
That is, all arrows in 2-cycles and loops are virtual and not part of the 
Gabriel quiver.

\subsection{Algebras with Gabriel quiver $3\cA$}

We take the quiver $Q$ obtained by glueing the second and the third type above,
identifying vertex $2$ with vertex $2'$.

\[
  \xymatrix{
    1
    \ar@(ld,ul)^{\ve}[]
    \ar@<.5ex>[r]^{\alpha}
    & 2
    \ar@<.5ex>[l]^{\beta}
    \ar@<.5ex>[r]^{\gamma}
    & 3
    \ar@<.5ex>[l]^{\delta}
    \ar@(ru,dr)^{\ve'}[]
  }
\]
The  permutation $g$ is of the form 
$(\alpha \ \gamma \ \delta \ \beta)(\ve)(\ve')$. \ 
Let $m_{\alpha} = k\geq 2$.\\
The case  $k=1$ is special, this gives the  {\it triangular} algebra,  called $T(\lambda)$ in \cite{WSA-GV}, here $\lambda \neq 1$.
With suitable choice of $c_{\bullet}$, the presentation of
the weighted surface algebra induces the 
 (Gabriel) presentation of $T(\lambda)$
 
\begin{align*}
 \alpha \beta\alpha  &= \alpha \gamma \delta ,&
  \delta \beta \alpha &= \lambda (\delta\gamma\delta)  ,&
 \beta\alpha\beta  &= \gamma\delta\beta ,
 &
 \beta\alpha\gamma  &=  \lambda(\gamma\delta\gamma),
  \\
 \alpha \beta \alpha \gamma  &= 0 ,&
 \beta \alpha\beta\alpha\beta &= 0 ,&
  \delta \gamma \delta \beta &= 0,&
  \gamma \delta \gamma \delta \gamma &= 0 ,
\\
  \alpha\beta \alpha \beta \alpha  &= 0 ,
 &
  \delta \gamma \delta \gamma \delta &= 0 ,&
  \delta \beta  \alpha \beta &= 0 .
\end{align*}
 
One can show that $T(\lambda)$ and $T(\mu)$ are not isomorphic for $\lambda\neq \mu$. 
 The weighted surface algebras with the same quiver and $\ve, \ve'$ virtual loops have Gabriel quiver
 denoted by  
  $Q(3\cA)_2^k$ in \cite{H} (which is the algebra with parameter $k$ in the the family 
 $Q(3\cA)_2$ of \cite{E}).

\subsection{Spherical algebras}

We have the algebras  whose quiver is  given by  the first building block where we identify $a_1=a_{n+1}$, for $n\geq 2$.
The case $n=2$ gives the algebra $S(\lambda)$, called spherical algebra,  introduced in \cite{WSA-GV}, Example 3.6, as follows.

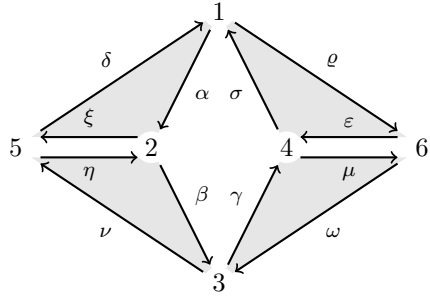
\[
\begin{tikzpicture}
[->,scale=.9]
\coordinate (1) at (0,2);
\coordinate (2) at (-1,0);
\coordinate (2u) at (-.925,.15);
\coordinate (2d) at (-.925,-.15);
\coordinate (3) at (0,-2);
\coordinate (4) at (1,0);
\coordinate (4u) at (.925,.15);
\coordinate (4d) at (.925,-.15);
\coordinate (5) at (-3,0);
\coordinate (5u) at (-2.775,.15);
\coordinate (5d) at (-2.775,-.15);
\coordinate (6) at (3,0);
\coordinate (6u) at (2.775,.15);
\coordinate (6d) at (2.775,-.15);
\fill[fill=gray!20] (1) -- (5u) -- (2u) -- cycle;
\fill[fill=gray!20] (1) -- (4u) -- (6u) -- cycle;
\fill[fill=gray!20] (2d) -- (5d) -- (3) -- cycle;
\fill[fill=gray!20] (3) -- (6d) -- (4d) -- cycle;
\node [fill=white,circle,minimum size=4.5] (1) at (0,2) {\ \quad};
\node [fill=white,circle,minimum size=4.5] (2) at (-1,0) {\ \quad};
\node [fill=white,circle,minimum size=4.5] (3) at (0,-2) {\ \quad};
\node [fill=white,circle,minimum size=4.5] (4) at (1,0) {\ \quad};
\node [fill=white,circle,minimum size=4.5] (5) at (-3,0) {\ \quad};
\node [fill=white,circle,minimum size=4.5] (6) at (3,0) {\ \quad};
\node (1) at (0,2) {1};
\node (2) at (-1,0) {2};
\node (2u) at (-1,0.15) {\ \quad};
\node (2d) at (-1,-0.15) {\ \quad};
\node (3) at (0,-2) {3};
\node (4) at (1,0) {4};
\node (4u) at (1,0.15) {\ \quad};
\node (4d) at (1,-0.15) {\ \quad};
\node (5) at (-3,0) {5};
\node (5u) at (-2.775,0.15) {};
\node (5d) at (-2.775,-0.15) {};
\node (6) at (3,0) {6};
\node (6u) at (2.775,0.15) {};
\node (6d) at (2.775,-0.15) {};
\draw[thick,->]
(1) edge node[below right]{\footnotesize$\alpha$} (2)
(2u) edge node[above]{\footnotesize$\xi$} (5u)
(5u) edge node[above left]{\footnotesize$\delta$} (1)
(5d) edge node[below]{\footnotesize$\eta$} (2d)
(2) edge node[above right]{\footnotesize$\beta$} (3)
(3) edge node[below left]{\footnotesize$\nu$} (5d)
(1) edge node[above right]{\footnotesize$\varrho$} (6u)
(6u) edge node[above]{\footnotesize$\varepsilon$} (4u)
(4) edge node[below left]{\footnotesize$\sigma$} (1)
(4d) edge node[below]{\footnotesize$\mu$} (6d)
(6d) edge node[below right]{\footnotesize$\omega$} (3)
(3) edge node[above left]{\footnotesize$\gamma$} (4)
;
\end{tikzpicture}
\]
where
the four shaded triangles denote the
$f$-orbits.
We take all multiplicities equal to $1$, the presentation induced by the
weighted surface algebra presentation is, with suitable choice of
$c_{\bullet}$, 

\begin{align*}
  \alpha \beta \nu &= \varrho \omega \nu ,
  &
  \beta \nu \delta &= \lambda \beta \gamma \sigma ,
  &
  \nu \delta \alpha &= \lambda \gamma \sigma \alpha ,
  &
  \delta \alpha \beta &=  \delta \varrho \omega ,
\\
  \gamma \sigma \varrho &=  \nu \delta \varrho ,
&
  \sigma \varrho \omega &= \lambda \sigma \alpha \beta ,
  &
 \varrho \omega \gamma &= \lambda \alpha \beta \gamma ,
  &
  \omega \gamma \sigma &= \omega \nu \delta ,
\\
  \alpha \beta \nu \delta \alpha &= 0 ,
  &
  \beta \nu \delta \varrho &= 0 ,
  &
  \nu \delta \alpha \beta \nu &= 0 ,
  &
  \delta \alpha \beta \gamma &= 0 ,
\\
 \gamma \sigma \varrho \omega \gamma &= 0 ,
  &
  \sigma \varrho \omega \nu &= 0 ,
  &
  \varrho \omega \gamma \sigma \varrho &= 0 ,
  &
  \omega \gamma \sigma \alpha &= 0 ,
\\
  \beta \gamma \sigma \varrho &= 0 ,
  &
  \sigma \alpha \beta \nu &= 0 ,
  &
  \delta \varrho \omega \gamma &= 0 ,
  &
  \omega \nu \delta \alpha &= 0 ,
\\
  \beta \nu \delta \alpha \beta &= 0 ,
  &
  \delta \alpha \beta \nu \delta &= 0 ,
 &
  \sigma \varrho \omega \gamma \sigma &= 0 ,
  &
  \omega \gamma \sigma \varrho \omega &= 0 .
\end{align*}

\subsection{The $n$-spherical algebra} When $n\geq 2$, the 
 permutation $g$ is of the form
$$ \prod_{i=1}^n (\xi_i \, \eta_i)\cdot
( \gamma_1\, \sigma_1\, \gamma_2\, \sigma_2\, \dots \, \gamma_n\, \sigma_n )\cdot
   ( \varrho_n\, \delta_n\, \varrho_{n-1}\, \delta_{n-1}\, \dots \, \varrho_1\, \delta_1 ).
$$
We take the multiplicities for  the $2n$-cycles to be $m, m'\geq 1$, and write $c, c'$ for the parameters
at these cycles.

\subsection{A mixed algebra} 

We can glue together the three building blocks by identifying 
 $2=a_1$, and $2'=a_{n+1}$.  
In this case, the permutation $g$ is the product of one large cycle with $n$ cycles of length $2$, and two loops:
$$\prod_{i=1}^n (\xi_i \ \eta_i)\cdot 
( \gamma_1\, \sigma_1\, \gamma_2\, \sigma_2\, \dots \, \gamma_n\, \sigma_n \ \gamma \ \delta \
    \varrho_n\, \delta_n\, \varrho_{n-1}\, \delta_{n-1}\, \dots \, \varrho_1\, \delta_1 \ \beta \ \alpha ) (\ve)(\ve').
$$
We take  again the multiplicities equal to $1$ on 2-cycles of $g$, or $m_{\gamma_1}=m$, and 
the parameter function with value $1$ on each virtual arrow. 
These algebras were not studied in previous papers but they fit into the same scheme.

\bigskip

\section{Construction of the module $M$}

Let $\La$ be one of the  algebras  as described above. 
Let $\Gamma\subset Q_0$ be the set of vertices which are not 
adjacent to a virtual arrow.

\begin{definition}
Let $M$ be the 
(right) $\La$-module
$$M:= \La \oplus [\bigoplus_{i\in \Gamma} S_i] \oplus [\bigoplus_{\nu\not\in \Gamma} \Omega^2(S_{\nu})]$$
\end{definition}
In the following we write down the details for the case of the $n$-spherical algebra, for the other algebras they are essentially the same.
In this case $\Gamma = \{ a_i\mid 1\leq i\leq n\}$. 

\bigskip

\subsection{The $\Omega$-translates of the simple modules}

For the algebra in question, the dimensions of the indecomposable projectives are:
$$\dim P_{a_i} = 2n(m+m'), \ \ \dim P_{b_i} = 2nm+2, \ \ \dim P_{d_i} = 2nm'+2.
$$
Let $a_i\in \Gamma$. The structure of $\Omega^{\pm 1}(S_{a_i})$ can be seen 
 from the presentation of the algebra.
The module $\Omega^2(S_{a_i})$  has dimension 5, the Loewy structure is
$$\begin{matrix} b_{i-1} && d_i \cr
         & a_i&\cr
b_i&& d_{i-1} \end{matrix}
$$
That is, the module has a 'simple waist'.
Now let  $\nu \in Q\setminus \Gamma$ we set $U_{\nu}:= \Omega^2(S_{\nu})$. 
Then $\Omega(U_{\nu})= \Omega^{-1}(S_{\nu})$ and $\Omega^{-1}(U_{\nu}) = \Omega(S_{\nu})$, their structure can also be seen  from the presentation.
We describe  $U_{\nu}$. 

\bigskip

\begin{lemma} The module $U_{b_i}$ is uniserial of length $2nm'-1$, with composition series  
$$U_{b_i} = \cU(a_i, d_{i-1}, a_{i-1}, d_{i-2}, a_{i-2}, \ldots, a_{i+1})
$$
The module $U_{d_i}$ is uniserial of length   $2nm-1$, with composition series  
$$ U_{d_i} = \cU(a_{i+1}, b_{i+1}, a_{i+2}, b_{i+2}, \ldots, a_i)
$$
(taking indices modulo $n$ and writing $a_i$, $d_i$, $b_i$ meaning the corresponding simple module).  
\end{lemma}

\bigskip

{\it Proof } We compute $U_{b_1}$, that is $\Omega^2(S_{b_1}) = \Omega(S_{b_1}) = \{ x\in P_{a_2} \mid \sigma_1x=0\} \subset P_{a_2}$. 
From the relations for the algebra, we have
$$\sigma_1\varrho_1\delta_1 = cA_{\sigma_1} = c\sigma_1A_{\sigma_1}'$$
Hence $\sigma_1\psi =0$ if we set
$$\psi= \psi_{\varrho_1\delta_1} := \varrho_1\delta_1 - cA_{\sigma_1}'.
$$
One exhibits a basis for  $\psi\La$, showing that it has  the same dimension as  $\Omega^2(S_{b_1})$, hence we have equality.
The submodule structure follows directly.   The case $U_{d_i}$ is similar.
$\Box$

\bigskip

\bigskip

\begin{proposition} We have ${\rm Ext}^1(M, M)=0$ and ${\rm Ext}^2(M, M)=0$.
\end{proposition}

{\it Proof } By ext symmetry, it suffices to show that for any non-projective indecomposable summand
$X$ of $M$ we have ${\rm Ext}^1(M, X)=0$ and ${\rm Ext}^1(X, M)=0$. 
For this, we use the following short exact sequences:

Let $a_i\in \Gamma$, 
$$0\to \Omega(S_{a_i})\to P_{a_i} \to S_{a_i}\to 0, \ \ 0\to \Omega^2(S_{a_i}) \ \to P_{c_i}\oplus P_{d_i} \to  \Omega(S_{a_i})\to 0
\leqno{(1)}
$$
Consider a vertex $\nu$ not in $\Gamma$, let $\nu=c_i$
$$0 \to \Omega^{-1}S_{b_i} \to P_{a_i} \to U_{b_i} \to 0, \ \ 0\to S_{b_i}\to P_{b_i} \to \Omega^{-1}(S_{b_i})\to 0
\leqno{(2)}
$$
Let $\nu=d_i$, then 
$$0 \to \Omega^{-1}(S_{d_i}) \to P_{a_i+1}  \to U_{d_i} \to 0, \ \ 0\to S_{d_i} \to P_{d_i} \to \Omega^{-1}(S_{d_i})\to 0
\leqno{(3)}
$$
We apply the functor ${\rm Hom}_A(-, X)$ to the above exact sequences. 

(I) Assume $X= S_{a_j}$ for some $j$. We know from the quiver that ${\rm Ext}^1(S_{a_i}, S_{a_j})=0$ already. To show that
${\rm Ext}^2(S_{a_i}, S_{a_j})=0$ we apply the functor to the second sequence in (1). From the structure of $\Omega^2(S_{a_i})$ we see directly that
${\rm Hom}(\Omega^2(S_{a_i}), S_{a_j})=0$ and hence ${\rm Ext}^2(S_{a_i}, S_{a_j})=0$. 

We have ${\rm Ext}^1(U_{b_i}, X)=0$
since ${\rm Hom}(\Omega^{-1}(S_{b_i}), S_{a_j})=0$ 
Furthermore ${\rm Ext}^2(U_{b_i}, X)=0$ since ${\rm Hom}(S_{b_i}, S_{a_j}) = 0$. 
Similarly one shows ${\rm Ext}^1(U_{d_i}, X)=0$ and ${\rm Ext}^2(U_{d_i}, X)=0$.

\medskip

(II) Now assume $X=U_{b_j}$ for some $j$. First, by dimension shift ${\rm Ext}^1(U_{\nu}, U_{b_j}) \cong {\rm Ext}^1(S_{\nu}, S_{b_j})=0$ 
for any $\nu$ of valency $1$, from the quiver. Next, consider
${\rm Ext}^2(U_{\nu}, X)$, by applying the functor ${\rm Hom}(-, X)$ to the second exact sequence in (2). 
We have ${\rm Hom}(S_{\mu}, U_{b_j})=0$ (the socle of $U_{b_i}$ is always some $S_a$), and hence 
${\rm Ext}^2(U_{\nu}, X)=0$. 

\medskip
Now consider ${\rm Ext}^t(S_{a_i}, X)$ for $t=1, 2$. By the ext symmetry, it is isomorphic to
${\rm Ext}^t(X, S_{a_i})$ for $t=2, 1$. By part (I) we know that it is zero.

The proof for $X=U_{d_i}$  is analogous.
$\Box$

\bigskip

\begin{remark}\normalfont  For possible later use, we write down sequences which may be used to show ${\rm Ext}^1(X, M)=0$ and ${\rm Ext}^2(X, M)=0$: 
Let $a_i\in \Gamma$, 
$$0\to S_{a_i}\to P_{a_i} \to \Omega^{-1}(S_{a_i}) \to 0, \ \ 0\to \Omega^{-1}(S_{a_i})\to P_{b_i}\oplus P_{d_i} \to  \Omega^{-2}(S_{a_i})\to 0
\leqno{(1^*)}
$$
Consider a vertex $\nu$ not in $\Gamma$, let $\nu=c_i$
$$0 \to  U_{b_i}  \to P_{a_{i+1}} \to \Omega(S_{b_i}) \to 0, \ \ 0\to \Omega(S_{b_i}) \to P_{b_i} \to S_{b_i}\to 0
\leqno{(2^*)}
$$
Let $\nu=d_i$, then 
$$0 \to  U_{d_i} \to   P_{a_{i}}  \to \Omega(S_{d_i}) \to 0, \ \   0 \to \Omega(S_{d_i}) \to P_{d_i} \to S_{d_i}\to 0.
\leqno{(3^*)}
$$
\end{remark}

\bigskip

\section{Ext vanishing and  3-cluster tilting}

We would like to determine when  $M$ is 3-cluster tilting. Hence take $X$ indecomposable and not projective, and assume
$${\rm Ext}^1(M, X)=0 = {\rm Ext}^2(M, X).
$$
By   ext symmetry, we get for free that
${\rm Ext}^1(X, M) = 0 = {\rm Ext}^2(X, M).$
The aim is to show that $X$ is in add$(M)$, or if not, to identify $X$.

\begin{lemma} 
The socle and the top of $X$ belong to ${\rm add}(\oplus_{i} S_{a_i})$.
\end{lemma}

{\it Proof} \  Let $\nu$ be a vertex $\neq a_i$ for any $i$.  Apply the functor ${\rm Hom}(-, X)$ to the second sequence of (2), this gives the exact sequence
$$0 \to {\rm Hom}(\Omega^{-1}(S_{b_i}), X ) \to {\rm Hom}(P_{b_i}, X) \to {\rm Hom} (S_{b_i}, X)\to 0
$$
Any homomorphsim $P_{b_i}\to X$ must map the socle to zero, otherwise
it would be split. Hence it lies in ${\rm Hom}(\Omega^{-1}(S_{b_i}), X)$
and therefore the  first two terms are isomorphic.
Hence the last term is zero, as required.
To show that also ${\rm Hom}(X, S_{\nu})=0$ we use a sequence from $(2^*)$. 
$\Box$

\bigskip

\begin{lemma}  We have ${\rm Hom}(\Omega(X), S_{a_i})=0$ and ${\rm Hom}(S_{a_i}, \Omega^{-1}(X))=0$.
\end{lemma}

\bigskip

{\it Proof } Since ${\rm Ext}^1(X, S_{a_i})=0$, from a minimal projective cover of $X$ we obtain the exact sequence
$$0\to {\rm Hom}(X, S_{a_i}) \to {\rm Hom}(P_X, S_{a_i}) \to {\rm Hom}(\Omega(X), S_{a_i})\to 0
$$
The first two terms are isomorphic since we start with a projective cover. Hence the last term is zero.
Similarly by using an injective hull we get ${\rm Hom}(\Omega^{-1}(X), S_{a_i}) = 0$.
$\Box$

\bigskip

Let $\cX$ be the category of $A$-modules which have socle and top in ${\rm add}(S_{a_i})$. This category is equivalent to
${\rm mod} -e\La e$. where  $e$ is the idempotent $e:= \sum_{i} e_{a_i}$. An equivalence is given by 
the functor $V\mapsto Ve$, with inverse the composite of $(-)\otimes_{e\La e}(e\La)$ follows by factoring out the largest $A$-submodule $V'$ with $V'e=0$
(see for example \cite{BDK}).

We may write down quiver and presentation of the algebra $e\La e$. The arrows are  $x_i := \gamma_i\sigma_i$ and 
$y_i:= \varrho_i\delta_i$, for $1\leq i\leq n$ where $x_i: a_i\mapsto a_{i+1}$ and $y_i: a_{i+1}\to a_i$. 
From the relations for $\La$ we see 
We claim that $x_iy_i=0$ and $y_ix_{i-1}=0$. That is,  $e\La e$ is special biserial. Moreover, for any $i$, the longest non-zero monomial 
$x_ix_{i+1}\ldots$ is up to a scalar equal to the longest non-zero monomial
$ y_{i-1}y_{i-2}\ldots$, and this gives the socle relations.

\begin{lemma} The module $X$ has simple socle and top.
\end{lemma}

{\it Proof }  The module $X$, and as well, all projectives (injectives) $P_{a_i}$ belong to the category $\cX$, and hence 
we may fix an injective hull, or projective cover, of $X$ by identifying with the image of a suitable
injective hull, or projective cover, of $Xe$,  in ${\rm mod} e\La e$. 
The indecomposable $e\La e$-modules are  'strings' or 'bands', and their injective hulls or projective covers may be written
down explicitly.

Assume the socle of $X$ is not simple, then consider the injective hull $I_X$, it has at least two indecomposable summands, say it
is $\oplus_{i\in R} P_{a_i}$. 
We may assume, with the above convention,
and taking $X\to I_X$ as inclusion, that $X$ has a generator $\omega = (\omega_1, \omega_2, 0, \ldots)$ such 
that $\omega A$ has socle of length two, and moreover, that $\omega x_j= (\omega_1x_j, 0, \ldots)$ and $\omega y_{j-1}= (0, \omega_2y_{j-1}, 0, \ldots)$ and
$\omega x_r=0$, $\omega y_s=0$ for all other generators $x_r, y_s$ of $e\La e$. This implies then that $\omega J\subseteq X$ where $J$ is the radical of $\La$. 
Now consider 
$\pi: I_X\to \Omega(X)$. The element  $\pi(\omega_1, 0, \ldots)$ is non-zero (since $\omega$ is a generator for $X$).
Furthermore, $[\pi(\omega_1, 0, \ldots)]J = \pi[(\omega_1, 0)J] =0$ since $(\omega_1, 0, \ldots ) J$ is contained in $X$.
Now $\pi(\omega_1, 0) = \pi(\omega_1, 0)e$,  (since ${\rm top } X$ is in add$(\oplus S_{a_i})$. Hence 
for some $i$ we have ${\rm Hom}(S_{a_i}, \Omega^{-1}(X))\neq 0$. This contradicts the previous Lemma.

Similarly by exploiting a projective cover, one shows that the top of $X$ must be simple.
$\Box$

\begin{proposition}  The module $X$ is uniserial.
\end{proposition}

\bigskip

{\it Proof } If $X$ is not uniserial then $Xe$ 
is not uniserial (using the structure of the projectives in this case). Then $Xe$ is a 'band module'. This means that
$X$ contains a submodule isomorphic to the second socle of some $P_{a_j}$. That is ${\rm Hom}(\Omega^2(S_{a_j}), X)\neq 0$. 

Applying  ${\rm Hom}(-, X)$ to the exact sequence
$$0\to \Omega^2(S_{a_j}) \stackrel{\iota}\to P:= P_{b_j}\oplus P_{d_{j-1}} \longrightarrow \Omega(S_{a_j}) \to 0
$$
gives an excact sequence, that is a non-zero homomorphism $\theta: \Omega^2(S_{a_j})\to X$ factors through 
$\iota$, say $\theta = \psi\circ \iota$.
The kernel of $\theta$ is the socle of $\Omega^2(S_{a_j})$ which also is the socle of $P$.  We  factor out these socles,  then for the induced maps we have
$$\bar{\theta} =  \bar{\psi}\circ \bar{\iota}.$$
Now, the map $\bar{\psi}$  on the socle of $\bar{P}$ is non-zero on each component. It follows that
the image of $\bar{\psi}$ has Loewy length equal to the Loewy length of $P/{\rm soc} P$. 

Note that all modules $P_{a_i}$ have the same Loewy length $\ell$ say. As well $P_{b_j}\oplus P_{d_{j-1}}$ has Loewy length $\ell$. 
Hence the Loewy length of $P/{\rm soc} P$ is $\ell-1$. 
The image of $\bar{\psi}$ is contained in the radical of $X$, which is the unique maximal submodule.
It follows that the Loewy length of $X$ is $\ell$. But this means that $X$ must be projective, a
contradiction.
This shows that $Xe$ is uniserial, and then from the structure of the projectives, also $X$ is uniserial.
$\Box$

\bigskip

We summarize. We have shown that if $X$ is indecomposable and not projective such that
${\rm Ext}^1(M, X)=0={\rm Ext}^2(M, X)$ then 

\medskip

$(*)$ \ $X$ is uniserial, and ${\rm soc} X$ and ${\rm top} X$ are in add$(\oplus_{i} S_{a_i})$. 
That is, $X$ is a subquotient of some $U_{b_i}$ or $U_{d_j}$. 

\bigskip

We show now that if $X$ is any module satifying $(*)$ 
then  ${\rm Ext}^1(M, X)=0$ and ${\rm Ext}^2(M, X)=0$. 

\begin{lemma} 
Let $X= \cU(a_j, b_j, a_{j+1}, b_{j+1}, \ldots a_l)$, a subquotient of some $U_{\nu}$. 
Then \\
 ${\rm Ext}^1(M, X) =0$ and ${\rm Ext}^2(M, X)=0$.
  \end{lemma}

{\it Proof}\  We use the sequences in the proof of Proposition 4.3.
We apply the functor $(-,X):= {\rm Hom}(-, X)$ to the exact sequences in (1). We start with the second, this gives
$$0\to  (\Omega(S_{a_i}), X) \longrightarrow (P_{b_i}\oplus P_{d_{i-1}}, X) \longrightarrow (\Omega^2(S_{a_i}), X) \to {\rm Ext}^1(\Omega(S_{a_i}), X)\to 0
$$
We see that ${\rm Hom}(\Omega^2(S_{a_i}), X)=0$ ($X$ is uniserial). Hence the ext space is zero. 
Moreover, it follows that 
the first two terms are isomorphic, which we can use for the first sequence:  
$$0 \to (S_{a_i}, X) \longrightarrow Xe_{a_i} \longrightarrow Xe_{b_i}\oplus Xe_{d_i}  \to {\rm Ext}^1(S_{a_i}, X)\to 0
$$
In our case, $Xe_{d_i}=0$. Note that in the composition series we have length two subquotients $a_r, b_r$, except that for $l=r$ we have an extra copy of $a_l$.
Hence if $i=l$ then the first term is $K$, and $|Xe_{a_l}| = 1 + |Xe_{b_l}|$ and ext is zero. 
Suppose $i\neq \ell$, then the first term is zero and the second and third are isomorphic. Again ext is zero.

\bigskip

Next, we apply $(-, X)$ to the sequences in (3). 
Since $S_{d_i}$ does not occur in $X$,  the functor takes the second sequence to zero.
From the first sequence we get 
$$0 \to {\rm Hom}(U_{d_i}, X) \to Xe_{a_i+1} \to 0 \to {\rm Ext}^1(U_{d_i}, X) \to 0
$$ and the ext space is zero.\\
Now consider $(-,X)$ applied to sequences in (2).
The second sequence gives ${\rm Hom}(\Omega^{-1}(S_{b_i}), X) \cong Xe_{b_i}$ and ${\rm Ext}^1(\Omega^{-1}(S_{b_i}), X)=0$. Consider the first
sequence, this gives 
$$0\to {\rm Hom}(U_{b_i}, X) \to Xe_{a_i} \to Xe_{b_i} \to {\rm Ext}^1(U_{b_i}, X)\to 0
$$
If the top (ie $S_{a_i}$) of $U_{b_i}$ is not the same as the socle of $X$ then the hom space is zero and the second and third term are isomorphic, and ext is zero.
Supoose $i=l$, then the first term is $K$, and $|Xe_{a_i}| = 1 + |Xe_{b_i}|$ and again the ext space  is zero.
$\Box$

\bigskip

\begin{corollary} Assume $\La$ is the triangle algebra, or the spherical algebra. 
Then $M$ is 3-cluster tilting.
\end{corollary}

{\it Proof} \  For these algebras, all indecomposables satisfying $(*)$ are in add$(M)$.
$\Box$

\bigskip

Consider an $n$-spherical algebra for $n\geq 3$ and $m=m'=1$. Then the (finite) set of modules 
$X$ satifying $(*)$ contains  all modules 
of the form 
$$\cU(a_i, b_i, a_{i+1}), \ \ \cU(a_i, d_{i-1}, a_{i-1}).
$$
To have a 3-cluster tilting  module with $M$ as a summand, we would need to take
$\wt{M} = M\oplus \cV$ where $\cV$ is the direct sum of all modules satisfying $(*)$.  However, 
$\wt{M}$ has self-extensions. For example there is a non-split exact sequence
$$0\to \cU(a_2, b_2, a_3) \to S_{a_2}\oplus \cU(a_1, b_1, a_2, b_2, a_3) \to \cU(a_1, b_1, a_2)\to 0
$$
Hence $M$ cannot be extended to a 3-cluster tilting  module for the $n$-spherical algebra when $n\geq 3$.

\bigskip

We also consider the algebra with triangular quiver and $k\geq 2$. In this case the list of
 uniserial modules $X$ which are subquotients of $U_1$ and $U_3$
 contains the modules
 $$\cU(2, 3, 2),  \ \cU(2, 1, 2)
$$
Let $\wt{M} = M\oplus \cV$ where $\cV$ is the direct sum of all indecomposable modules satisfying $(*)$. 
This is not a 3-cluster tilting module since it has self-extensions:
we have the non-split exact sequence
$$0\to \cU(2, 3, 2)\to S_2\oplus \cU(2, 1, 2, 3, 2)\to \cU(2, 1, 2)\to 0
$$

 \end{document}